\documentclass[12pt]{article} 
\usepackage{amsmath, amssymb}
\usepackage{amsthm} 
\theoremstyle{plain} 
\newtheorem{theorem}{\indent\sc Theorem}[section]
\newtheorem{lemma}[theorem]{\indent\sc Lemma}
\newtheorem{corollary}[theorem]{\indent\sc Corollary}

\theoremstyle{definition} 
\newtheorem{definition}[theorem]{\indent\sc Definition}
\newtheorem{remark}[theorem]{\indent\sc Remark}

\makeatletter
\def\address#1#2{\begingroup
\noindent\parbox[t]{7.8cm}{%
\small{\scshape\ignorespaces#1}\par\vskip1ex
\noindent\small{\itshape E-mail address}%
\/: #2\par\vskip4ex}\hfill%
\endgroup}%
\makeatother
\title{\uppercase {$LS_{\lowercase {r}}$-valued Gauss maps\\ and \\ spacelike  surfaces of revolution in $\mathbb R_1^4$
}} 
\author{
\bigskip \\
\textsc{Dang Van Cuong}\footnote{ \ The author is supported in part by the National Foundation for Science and Technology
Development, Vietnam (Grant No. 101.01.30.09)}}
\date{\today}
\begin{document}
\maketitle

\footnote{AMS
\textit{Mathematics Subject Classification}.
53C50; 53A35.
}
\begin{abstract}
 To study  spacelike surfaces in the Lorentz-Minkowski space $\Bbb R^{4}_1,$ we construct a pair of maps, called $\mathfrak l_r^{\pm} $-Gauss maps, whose values are in the lightcone. We can use these maps to study umbilical spacelike surfaces and find parametrizations of spacelike surfaces of revolution of hyperbolic and elliptic types in some particular cases.
\end{abstract}
\section{Introduction}
 It is well known that  Gauss maps have been used as one of the most useful tools to study the behavior and geometric invariants of hypersurfaces. The Gauss maps associated with an arbitrary normal field  $\nu$ are  particularly important for the study of surfaces with codimension larger than one. This method provides a convenient way to study invariants and characters of surfaces that are dependent or independent on $\nu$.  For instance, Izumiya et al. \cite{izu1} used the Gauss map associated with a normal field
to study  $\nu$-umbilicity for spacelike surfaces of codimension two in the Lorentz-Minkowski space. Marek Kossowski \cite{ko}, on the other hand, introduced two $S^2$-valued Gauss maps, whose values are in the lightcone, in order to study spacelike 2-surfaces in $\Bbb R^4_1$, and the method was  continued by Izumiya et al. \cite{izu2} to study  properties of spacelike surfaces of codimension two.  Introducing the $HS_r$-valued Gauss maps, called $\textbf  n_r^{\pm}$-Gauss maps, Cuong and Hieu \cite{C-H} has studied  the umbilical spacelike surfaces of codimension two in $\mathbb R_1^{n+1}$.  By a similar way, it is possible to introduce  $LS_r$-valued Gauss maps, called $\mathfrak l_r^{\pm}$ -Gauss maps, for the study of spacelike surfaces of codimension two.\\
\indent As an application of the  $\mathfrak l_r^{\pm} $-Gauss maps, in this paper we are going to use $LS_r$-valued Gauss maps in order to study umbilcal spacelike surfaces of codimension two, especially when the surfaces lie in de Sitter.  Using $\mathfrak l_r^{\pm}$-Gauss maps, we are able to calculate the curvatures of spacelike surfaces of revolution of hyperbolic or elliptic types, which are the orbits of a regular spacelike curve under the orthogonal transformations of $\mathbb R_1^4$ leaving a spacelike or timelike plane. We then specifically describe the totally umbilical and maximal spacelike surfaces of revolution and give  explicitly their parametrizations. The proofs of these results are mainly based on two useful lemmas that involve solving differential equations. From now on we assume that the  surface $M$ is a spacelike surface of codimension two.

\section{Preliminary Notes}
\subsection{The Lorentz-Minkowski Space $\Bbb R^{4}_1$}
 The Lorentz-Minkowski space $\Bbb R^{4}_1$ is the $4$-dimensional vector space $$\Bbb R^{4}=\{( x_1, x_2,x_3,x_4): x_i\in \Bbb R, i=1,\ldots, 4\},$$ endowed with the pseudo scalar product  defined by
$$\langle \textbf x, \textbf y\rangle=\sum_{i=1}^{3}x_iy_i-x_{4}y_{4},$$
where $\textbf x=( x_1,\ldots, x_{4}), \textbf y=(y_1, \ldots y_{4})\in \Bbb R^{4}.$
Since $\langle \cdot,\cdot \rangle$ is non-positively defined, $\langle \textbf x, \textbf x\rangle$ may receive negative value, , for a given $\textbf x\in \Bbb R^4_1$. We say that a nonzero vector $\textbf x\in \Bbb R^{4}_1$ is spacelike, lightlike or timelike if  $\langle \textbf x, \textbf x\rangle>0$, $\langle \textbf x, \textbf x\rangle=0$ or $\langle \textbf x, \textbf x\rangle<0$, respectively.
 Two  vectors $\textbf  x, \textbf  y \in \Bbb R^4_1$ are said to be pseudo-orthogonal if  $\langle \textbf x, \textbf y\rangle=0$.

The norm of a vector $\textbf x\in \Bbb R^{4}_1$ is defined as follows $$\|\textbf x\| =\sqrt{|\langle \textbf x,\textbf x\rangle|}.$$

For a nonzero vector $\textbf n\in \Bbb R^{4}_1$ and a contant $c\in \Bbb R$, the hyperplane with pseudo normal $\textbf n$ is defined by
$$HP(\textbf n,c)=\{\textbf x\in\Bbb R^{4}_1 : \langle \textbf x,\textbf n\rangle=c\}.$$

The hyperplane $HP(\textbf n,c)$ is called spacelike, lightlike or timelike if $\textbf n$ is timelike, lightlike or spacelike, respectively.

For a vector $\textbf a\in\Bbb R^4_1$ and a positive constant $R$, we define:
\begin{enumerate}
\item The Hyperbolic with center $\textbf a$ and radius $R$ by
  $$H^{3}(\textbf a,R)=\{\textbf x\in\Bbb R^{4}_1\ |\ \langle \textbf x-\textbf a,\textbf x-\textbf a\rangle=-R\}.$$
\item The de Sitter with center $\textbf a$ and radius  $R$ by
$$S_1^{3}(\textbf a,R)=\{\textbf x\in\Bbb R^{4}_1\ |\ \langle \textbf x-\textbf a,\textbf x-\textbf a\rangle=R\}.$$
\item The Lightcone with vertex $\textbf  a$ by
$$LC^3(\textbf a)=\{\textbf x\in\Bbb R^{4}_1:  \langle \textbf x-\textbf a,\textbf x-\textbf a\rangle=0\},$$
if $\textbf  a=0$ we have Lightcone $LC^3$.
\end{enumerate}

\subsection{$\mathfrak l_r^{\pm}$-Gauss Maps}
 In this section, we will introduce the Gauss maps of a surface and its curvature.

Let $M$ be a spacelike surface with parametrization
  $\textbf X:U\to\Bbb R^{4}_1$ , where $U\subset\Bbb R^{2}$ is open and  $M=\textbf X(U)$ is identified with $\textbf X.$  The normal plane of $M$ at $p\in M,$ denoted by $N_pM,$  can be viewed as a timelike 2-plane passing the origin. It is easy to see that $N_pM$ intersects with the lightcone $LC^3$ by two lines. Fix $r>0$, the hyperplane $\{x_{4}=r\}$ meets these lines exactly at two points, denoted by $\mathfrak{l}_r^{\pm}(p).$\\
\indent  The maps $$\begin{aligned}\mathfrak{l}^{\pm}_r:\  &M&&\to&& LS_r:=LC^3\cap \{x_{4}=r\},&\\ &p&&\mapsto& &\mathfrak{l}^{\pm}_r(p)&\end{aligned}$$ are called  $\mathfrak l^{\pm}_r$-Gauss maps.\\
\indent Algebraically, $\mathfrak l_r^{\pm} $ are the solutions of  the following system of  equations
\begin{equation} \label{hel} \left\{\begin{aligned}
&\langle \mathfrak l,\textbf  X_{u_i}\rangle =0,\ i=1,2\\
& \langle \mathfrak l,\mathfrak l\rangle =0,\\
&\mathfrak l_4=r.
\end{aligned} \right.
 \end{equation}
\indent It is easy to show that  $\mathfrak l_r^{\pm} $ are smooth vector fields, see \cite{C-H}. Therefore,  using  the symbol `` * '' instead of  `` + '' or `` - '',   we have the following notions and facts.\\
\indent The derivative of  $\mathfrak l_r^{*} $ at $p,$
$$d\mathfrak l_r^{*} (p)\ :\ T_pM\rightarrow T_{\mathfrak l_r^{*} (p)}LC^3\subset T_pM\oplus N_pM,$$
 can be  descomposed as
$$d\mathfrak l_r^{*}(p)=d{\mathfrak l_r^{*}}^T(p)+d{\mathfrak l_r^{*}}^N(p),$$
where
$d{\mathfrak l_r^{*}}^T$ and $d{\mathfrak l_r^{*}}^N$ are the tangent and normal components of $d{\mathfrak l_r^{*}},$ respectively.\\
\indent So,  we define:
\begin{enumerate}
\item  The {\it $\mathfrak l_r^{*}$-Weingarten map} of $M$  at $p$ by
$$A_p^{\mathfrak l_r^{*}}:=-d{\mathfrak l_r^{*}}^T(p).$$
\item The $\mathfrak l_r^{*}$-Gauss-Kronecker curvature of $M$  at $p$ by
$$K_p^{\mathfrak l_r^{*}}:=\det(A_p^{\mathfrak l_r^{*}}).$$
\item The $\mathfrak l_r^{*}$-mean curvature  of $M$  at $p$ by
$$H^{\mathfrak l_r^{*}}_p:=\frac{1}{2}\text{trace}(A_p^{\mathfrak l_r^{*}}).$$
\item The $\mathfrak l_r^{*}$-principal curvatures of $M$ at $p$ by the eigenvalues   of $A_p^{\mathfrak l_r^{*}}$,
$$k_1^{\mathfrak l_r^{*}}(p),k_2^{\mathfrak l_r^{*}}(p).$$
\end{enumerate}
It is clear that
$$K_p^{\mathfrak l_r^{*}}=k_1^{\mathfrak l_r^{*}}(p)k_{2}^{\mathfrak l_r^{*}}(p),$$
and
$$H_p^{\mathfrak l_r^{*}}=\frac 1{2}(k_1^{\mathfrak l_r^{*}}(p)+k_{2}^{\mathfrak l_r^{*}}(p)).$$

 Furthermore, we have some well-known facts:
\begin{enumerate}
\item  The $\mathfrak l_r^{*}$-Weingarten map is self-adjoint.
\item  The $\mathfrak l_r^{*}$-principal curvatures  $k_1^{\mathfrak l_r^{*}}(p)$ and $k_2^{\mathfrak l_r^{*}}(p)$ of $M$ at $p$  are the solutions of the following equation
\begin{equation} \label{principal} \det(b_{ij}^{\mathfrak l_r^{*}}(p)-kg_{ij}(p))=0,\end{equation}
where  $b_{ij}^{\mathfrak l_r^{*}}(p):=\langle \textbf X_{u_iu_j}(p),\mathfrak l_r^{*}(p)\rangle,\ i,j=1,2, $ are the coefficients of the $\mathfrak l_r^{*}$-second fundamental form of $M$ at $p$.
\item $K_p^{\mathfrak l_r^{*}}={\det(b_{ij}^{\mathfrak l_r^{*}}(p))}.{\det(g_{ij}(p))}^{-1}.$
\end{enumerate}

\indent Now, let $\{e_3,e_4\}$ be an orthonormal frame of the normal bundle of $M$, and let $A^{e_3}$,  $A^{e_4}$ be the shape operators associated to $e_3,e_4$, respectively, where $\langle e_3,e_3\rangle=1,\langle e_4,e_4\rangle=-1$. The mean curvature vector $H$ is given by
$$H=\frac{1}{2}(\text{trace}A^{e_3})e_3-\frac{1}{2}(\text{trace}A^{e_4})e_4=H^{e_3}e_3-H^{e_4}e_4.  $$
It is easily to see that $H=0$ iff  $H^{e_3}=0$ and $H^{e_4}=0$, which is equivalent to  $H^{\mathfrak l_r^{+}}=0$ and $H^{\mathfrak l_r^{-}}=0$.
\begin{definition}\
\begin{enumerate}
\item A point $p\in M$ is said to be $\mathfrak l_r^{*} $-umbilic if $k_1^{\mathfrak l_r^{*} }(p)=k_2^{\mathfrak l_r^{*} }(p)=k.$ In the case $k=0$ we said that $p$ is $\mathfrak l_r^{*} $-flat.
\item $M$ is said to be $\mathfrak l_r^{*} $-umbilic $($$\mathfrak l_r^{*}$-flat$)$ if each point $p\in M$ is $\mathfrak l_r^{*} $-umbilic $($$\mathfrak l_r^{*}$-flat$)$.
\item $M$ is said to be totally umbilical $($totally flat$)$ if each point $p\in M$ is $\mathfrak l_r^{*}$-umbilic $($$\mathfrak l_r^{*} $-flat$)$ for every $r>0.$
\item $M$ is said to be maximal if $H=0$.
 \end{enumerate}
\end{definition}
\section{Totally Umbilical Spacelike Surfaces}
In this section, we study the umbilicity of a surface by using the $\mathfrak l_r^{*} $-Gauss map with respect to a fixed $r$. \\
Following the approach developed in the proof of Theorem 3.2 in \cite{C-H}, we are able to obtain a similar result. Our method is based on the $\mathfrak l_r^{*} $-Gauss map on a lightcone instead of using the $\textbf  n_r^{*} $-Gauss map on a Hyperbolic.
\begin{theorem} Let $M$ be a connected spacelike surface. The following statements are equivalent:
\begin{enumerate}
\item There exists a constant $r>0$ such that $M$ is $\mathfrak l_r^{*} $-flat.
\item There exists a constant $r>0$ such that $\mathfrak l_r^{*} $ is constant.
\item There exist a lightlike vector $\text{\bf n}=(n_1,n_2,n_3,n_4)$ and a real number $c$ such that $M\subset HP(\text{\bf n},c)$.
\end{enumerate}
\end{theorem}
\indent Using $HS_r$-valued Gauss maps, Cuong and Hieu \cite{C-H} showed that a surface is totally umbilical iff it is $\nu$-umbilic for every smooth normal vector field $\nu$. Then the umbilicity of surface was studied by using $\textbf  n_r^{*}$-Gauss maps, especially when the surface lies in $H^n(0,R)$. Some of these results are based on the linear independence of the position vector $\textbf X$ and the normal vector field $\textbf  n_r^{*} $. Unfortunately, when the surface is contained in $S_1^3(0,R)$,  the pair of vectors $\textbf X$ and $\textbf  n_r^{*} $ may be linearly dependent.  Therefore, we use $\mathfrak l_r^{*}$ instead of $\textbf  n_r^{*}$ in order to study the surface in de Sitter. The reader is referred to  \cite{C-H} for more details of  the proof.
\begin{corollary}\label{inS} If $M$ is contained in the intersection of a de Sitter and a hyperplane, then $M$ is totally umbilic.
\end{corollary}
\begin{theorem}\label{inhyper}
Let $M$ be a surface  in $S_1^3(0,R)$. The following statements are equivalent:
\begin{enumerate}
\item There exists $r>0$ such that $M$ is $\mathfrak l_r^{*} $-umbilic.
\item $M$ is totally umbilic.
\item $M$ is contained in a hyperplane.
\end{enumerate}
\end{theorem}
Using $\mathfrak l_r^{*} $-Gauss maps, we can find some necessary and sufficient conditions for a surface lying in a de Sitter to be a part of a horizontal hyperplane (i.e. the hyperplane is pseudo orthogonal to the time-axis). The reader is referred to Theorem 4.7 in \cite{C-H} for more details of  the proof.
\begin{theorem}\label{hori} Let $M$ be a surface contained in $S_1^3(\textbf  a,R)$. The following statements are equivalent:
\begin{enumerate}
\item $M$ is contained in a hyperplane  $HP(\text{\bf  n},c)$, where $\text{\bf  n}=(0,0,0,1)$.
\item $\mathfrak l_r^{*}$ is parallel for any $r>0$.
\item There exist two different parallel normal fields $\mathfrak l_{r_1}^{*}, \mathfrak l_{r_2}^{*}$.
\item There exists $r>0$ such that $A^{\mathfrak l_r^{*}}=-\lambda id$, where $\lambda$ is a constant.
\end{enumerate}
\end{theorem}
The following theorem gives another necessary and sufficient conditions  for a surface to be a part of the intersection of $S_1^3(\textbf  a,R)$ and a horizontal hyperplane  without assuming that it lies in the de Sitter. The reader is referred to Theorem 4.8 in  \cite{C-H} for more details of  the proof.
\begin{theorem} Let $M$ be a  surface in $\Bbb R_1^4$. The following statements are equivalent:
\begin{enumerate}
\item There exists $r>0$ such that $\mathfrak l_r^{*} $ is parallel but not a constant vector, and $M$ is $\mathfrak l_r^{*} $-umbilic.
\item $M$ is contained in the intersection of $S_1^3(\textbf  a,R)$ and a horizontal hypeplane $HP(\text{\bf  n},c)$.
\end{enumerate}
\end{theorem}
\indent We will now  study totally umbilical surfaces  in the  general case. Let us recall some imporant notions and facts.\\
\indent We denote by  $\nabla^{\bot}$  the normal connection on the surface $M$. The normal curvature of $M$ is defined by
$$\begin{aligned}R^{\bot}:&\mathfrak X(M)\times\mathfrak X(M)\times \mathfrak X(M)^{\bot}&&\to&&\mathfrak X(M)^{\bot}&\\&\left(V,W,X\right)&&\mapsto&&\nabla_{[V,W]}^{\bot}X-\left[\nabla_V^{\bot},\nabla_W^{\bot}\right]X,&\end{aligned}$$
where $\mathfrak X(M)$ is the set of  smooth tangent vector fields and $ \mathfrak X(M)^{\bot}$ is the  set of smooth  normal vector fields on $M$, see \cite{os}.\\
\indent It can be shown, as a consequence of the Ricci equation \cite[p.125]{os}, that if $p$ is an umbilical point for some normal field $\nu$, then $R^{\bot}(p)=0$. Therefore, if $M$ is totally umbilical, then $R^{\bot}=0$ on $M$.
\begin{definition}{\rm \cite[p.6]{ct}} A connection is called flat if its curvature is zero.
\end{definition}
Proposition 1.1.5 in \cite{ct} gives us the following corollary.
\begin{corollary}\label{parallel} If $M$ is totally umbilic, then for any $p\in M$ there exist a local neighbourhood $U_p\subset M$ of $p$ and a parallel  normal frame, $\{\text{\bf  u},\text{\bf  v}\}$ on $U_p$.
\end{corollary}
Using the notion in the Corollary \ref{parallel}, we have the following result.
\begin{theorem} Let $M$ be a connected surface, then the following statements are equivalent:
\begin{enumerate}
\item $M$ is totally umbilic.
\item For any $p\in M$, $U_p$ is contained in the intersection of a hyperplane and either a Hyperbolic or a de Sitter.
\end{enumerate}
\end{theorem}
\begin{proof}\ \\
$(1)\Rightarrow (2)$.   Sine $\textbf  u$ is parallel, its norm is constant, therefore it is either lightlike, spacelike or timelike on $U_p$. Similarly to the vector $\textbf  v$.
\begin{enumerate}
\item[(i)] If $\textbf  u$ is timelike then $U_p$ is contained in the intersection of a hyperbolic and a hyperplane, by virtue Theorem 4.3 in \cite{izu1} and Theorem 4.5 in \cite{C-H}.
\item[(ii)] If $\textbf  u$ is spacelike then $U_p$ is contained in the intersection of a de Sitter and a hyperplane, by virtue Theorem 4.3 in \cite{izu1} and Theorem \ref{hori}.
\item[(iii)] Ne now consider the case when both $\textbf  u$ and $\textbf  v$ are lightlike. Set
$$\textbf  Z=\frac{\textbf  u}{\langle \textbf  u,\textbf  v\rangle }-\textbf  v.$$
Then $\textbf  Z$ is a smooth normal field on $U_p$ and
$$\langle \textbf  Z,\textbf Z\rangle=-2,\ \langle \textbf  Z,\textbf  v\rangle=1.  $$
Therefore,
$$\langle d\textbf  Z,\textbf  Z\rangle=0,\ \langle d\textbf  Z,\textbf  v\rangle=-\langle \textbf  Z,d\textbf  v\rangle= 0. $$
That means $\textbf  Z$ is parallel on $U_p$. Similarly to (i)  we imply that $U_p$ is contained in the intersection of a hyperbolic and a hyperplane.
\end{enumerate}
$(2)\Rightarrow (1).$  It follows from Theorem 4.5 in \cite{C-H} and Theorem \ref{hori}.
\end{proof}
\section{The Spacelike Surfaces of Revolution}
In this section we will use $\mathfrak l_r^{*} $-Gauss map to study totally umbilical and maximal  spacelike surfaces of revolution of hyperbolic and elliptic types in $\mathbb R_1^4$.
\subsection{Useful Lemmas}
We begin with two useful lemmas.
\begin{lemma}\label{hplane} If $f'\ne 0$ and $g''f'-g'f''=0$, then $g=cf+k$, where $c,k$ are constants.
\end{lemma}
\begin{proof} We have
$$\frac{g''f'-g'f''}{(f')^2}=0\ \Rightarrow \left(\frac{g'}{f'} \right)'=0\ \Rightarrow g=cf+k, $$
where $c,k$ are constants.
\end{proof}
\begin{lemma} \label{root}The solutions of the equation
$$\rho(u).\rho''(u)-(\rho'(u))^2-1=0$$
are
$$\rho(u)=\pm\frac{1}{\sqrt C} \cosh(\sqrt C(u-C_1)),$$
where $C>0,C_2$ are constants.
\end{lemma}
\begin{proof}
Set $t=\rho'(u)$. Since
$$\rho''(u)=\frac{(\rho'(u))^2}{\rho(u)}+\frac{1}{\rho(u)},$$
we have
$$t'(\rho)=(t(\rho))^{-1}\left(\frac{(t(\rho))^2}{\rho(u)}+\frac{1}{\rho(u)}  \right).$$
That means
\begin{equation}\label{ebn1}t'(\rho)-\frac{1}{\rho}t(\rho)=\frac{1}{\rho}(t(\rho))^{-1}.  \end{equation}
Equation (\ref{ebn1}) is a Bernoulli equation and its solutions are
$$t=\pm\sqrt{-1+C\rho^2},\ C=const>0.$$
That means
$$\rho'(u)=\pm\sqrt{-1+C(\rho(u))^2}.$$
The last equation gives us the result of the lemma.
\end{proof}
\begin{remark}\label{rm1}Similarly, we can show that:
\begin{enumerate}
\item The solutions of the equation
$$\rho(u).\rho''(u)+(\rho'(u))^2+1=0$$
are
$$\rho(u)=\pm\sqrt{C-(u-C_1)^2} ,$$
where $C_1,C_2$ are constants.
\item The solutions of the equation
$$\rho(u).\rho''(u)-(\rho'(u))^2+1=0$$
are
$$\rho(u)=\pm\frac{1}{\sqrt C} \sinh(\sqrt C(u-C_1)),  $$
where $C_1,C_2>0$ are constants.
\item The solutions of the equation
$$\rho(u).\rho''(u)+(\rho'(u))^2-1=0$$
are
$$\rho(u)=\pm\sqrt{(u-C_1)^2-C_2},$$
where $C_1,C_2$ are constants.
\end{enumerate}
\end{remark}
\subsection{Spacelike Surface of Revolution of Hyperbolic Type}
Let $C$ be a spacelike curve in $\text{span}\{e_1,e_2,e_4\}$ parametrized by  arc-length,
$$z(u)=\left(f(u),g(u),0,\rho(u)\right),\ u\in I.$$
The orbit of $C$ under the action of the orthogonal transformations of $\mathbb R_1^4$  leaving the spacelike plane $Oxy,$
$$A_S=\left[\begin{matrix}1&0&0&0\\0&1&0&0\\0&0&\cosh v&\sinh v\\0&0&\sinh v&\cosh v\end{matrix}\right],\ v\in\mathbb R,$$
is a  surface $M$ given by
\begin{equation}\label{hr1}
{\rm{\bf X}}(u,v)=\left(f(u),g(u),\rho(u)\sinh v, \rho(u)\cosh v\right),\ u\in I,\ v\in\mathbb R.
\end{equation}
The first partial derivatives  of ${\rm{\bf X}}(u,v)$ can be calculated by
$${\rm{\bf X}}_u=\left(f'(u),g'(u),\rho'(u)\sinh v,\rho'(u)\cosh v\right),\ {\rm{\bf X}}_v=\left(0,0,\rho(u)\cosh v, \rho(u)\sinh v\right),$$
and the coefficients of the first fundamental form of $M$ are
$$g_{11}=(f'(u))^2+(g'(u))^2-(\rho'(u))^2=1,\  g_{12}=0,\ g_{22}=\left(\rho(u)\right)^2>0. $$
It follows that $M$ is a spacelike surface, which is called  the {\it spacelike surface of revolution of hyperbolic type} in $\mathbb R_1^4$. From now on we always assume that $f'\ne0,g'\ne0$ and $\rho'\ne 0$.\\
The system of equations (\ref{hel}) yields $\mathfrak l_r^{\pm} =(\mathfrak l_1,\mathfrak l_2,\mathfrak l_3,r) $, where
$$\mathfrak l_1= \frac{r\rho'}{f'\cosh v}-\frac{g'}{f'}\mathfrak l_2,$$
$$\mathfrak l_2=\frac{g'\rho'\pm rf'}{\cosh v[(f')^2+(g')^2]},$$
$$\mathfrak l_3=r\tanh v.$$
Then  the  coefficients of the second fundamental form of $M$ associated to $\mathfrak l_r^{\pm} $-Gauss maps  are defined below:
$$b_{11}^{\mathfrak l_r^{\pm}  }=\langle \mathfrak l_r^{\pm}  ,{\rm{\bf X}}_{uu}\rangle=f'' \mathfrak l_1+g''\mathfrak l_2-\frac{r\rho'' }{\cosh v}, $$
$$b_{12}^{\mathfrak l_r^{\pm}}=\langle \mathfrak l_r^{\pm}  ,{\rm{\bf X}}_{uv}\rangle= r\rho'\cosh v.\tanh v-r\rho'\sinh v=0,$$
$$b_{22}^{\mathfrak l_r^{\pm} }=\langle \mathfrak l_r^{\pm} ,{\rm{\bf X}}_{vv} \rangle =-\frac{r\rho}{\cosh v}. $$
Solving the equation $\det \left[(b_{ij}^{\mathfrak l_r^{\pm}})-k(g_{ij})\right]=0$, we obtain the following principal curvatures of the surface in terms of the $\mathfrak l_r^{\pm} $-Gauss maps
\begin{equation}\label{k1+}k_1^{\mathfrak l_r^{+}}=\frac{rf''\rho'}{f'\cosh v}-r\frac{\rho'' }{\cosh v}+r\frac{g''f'-g'f''}{f'\left[(f')^2+(g')^2\right]\cosh v}\left(g'\rho' +f' \right), \end{equation}
\begin{equation}\label{k1-}k_1^{\mathfrak l_r^{-}  }=\frac{rf''\rho'}{f'\cosh v}-\frac{r\rho''}{\cosh v}+r\frac{g''f'-g'f''}{f'\left[(f')^2+(g')^2\right]\cosh v}\left(g'\rho' -f' \right),  \end{equation}
\begin{equation}\label{k2+-}k_2^{\mathfrak l_r^{+} }=k_2^{\mathfrak l_r^{-}}=-\frac{r}{\rho(\cosh v)}.\end{equation}
\begin{theorem}\label{u1} If the  surface defined by (\ref{hr1}) is totally umbilic, then it is contained in a timelike hyperplane and  parametrized by $$f(u)=\pm\frac{1}{\sqrt C}\sinh \left(\sqrt C(u-C_1)\right)+m ,$$
$$g(u)=\frac{C_2}{\sqrt C}\sinh \left(\sqrt C(u-C_1)\right) +k,  $$
$$\rho(u)=\pm\frac{1}{\sqrt C}\cosh \left(\sqrt C(u-C_1)\right) ,$$
where $C>0,C_1,C_2,k$ are constants.
\end{theorem}
\begin{proof}
From (\ref{k1+}), (\ref{k1-}) and (\ref{k2+-}), the umbilical condition for $M$ gives us the following equations:
$$k_1^{\mathfrak l_r^{+} }=k_2^{\mathfrak l_r^{+} }=k_2^{\mathfrak l_r^{-}}=k_1^{\mathfrak l_r^{-}}.$$
That mean we have a system of equations
$$\left\{\begin{aligned}&\frac{g''f'-g'f''}{f'}=0,\\
&\frac{f''\rho'-f'\rho''}{f'}=-\frac{1}{\rho},\\
&(f')^2+(g')^2-(\rho')^2=1.\end{aligned}\right. $$
The first equation tells us that  $M$ lies in a timelike hyperplane. By calculating $f$ and $g$ in terms of $\rho$, we have an equivalent system of equations
$$\left\{\begin{aligned}&g=C_2f+k\\&f'=\sqrt{\frac{1+(\rho')^2}{1+C^2} }\\
&\rho.\rho''-(\rho')^2-1=0.
\end{aligned}\right. $$
The conclusion of the theorem follows from  Lemma \ref{hplane} and Lemma \ref{root}.
\end{proof}
\begin{remark} By changing the parameters
$$t=\sqrt C(u-C_1),\ v=v,$$
we get a new parametrization of the totally umbilical spacelike  surface of revolution defined by (\ref{hr1})
$${\rm{\bf X}}(t,v)=\left(A\sinh t +m,B\sinh t +k,A\cosh t\sinh v, A\cosh t\cosh v\right) ,$$
where $A=\pm\frac{1}{\sqrt C} ,B=\frac{C_2}{\sqrt C} ,k,m $ are constants.
\end{remark}
\begin{theorem} If the  surface  defined by (\ref{hr1}) is  maximal, then it is contained in a timelike hyperplane and parametrized by
$$ f(u)=\frac{\sqrt{C_3}}{1+C_2^2}\arcsin\left(\frac{u-C_1}{\sqrt{C_3}} \right)+m,$$
$$g(u)=C_2f(u)+k,$$
$$\rho(u)=\pm\sqrt{C_3-(u-C_1)^2},$$
where $C_1,C_2,C_3>0,k$ are constants.
\end{theorem}
\begin{proof}
It follows from (\ref{k1+}), (\ref{k1-}) and (\ref{k2+-}) that the conditions for the maximality of $M$ are described by the following equations
$$k_1^{\mathfrak l_r^{+}}=-k_2^{\mathfrak l_r^{+} }=-k_2^{\mathfrak l_r^{-} }=k_1^{\mathfrak l_r^{-}}.$$
Therefore, we have the following system of equations
$$\left\{\begin{aligned}&\frac{g''f'-g'f''}{f'}=0,\\
&\frac{f''\rho'-f'\rho''}{f'}=\frac{1}{\rho},\\
&(f')^2+(g')^2-(\rho')^2=1.  \end{aligned}\right. $$
Similar to the proof of Theorem \ref{u1}, the conclusion of the theorem follows from  Remark \ref{rm1}.
\end{proof}
\begin{remark}
By changing the parameters
$$t=\sin\frac{u-C_1}{\sqrt C_3} ,\ v=v,$$
we get a new parametrization of the maximal spacelike surface of revolution defined by (\ref{hr1})
$${\rm{\bf X}}(t,v)=\left(A t+m,B t+k, C\cos t \sinh v, C\cos t \cosh v\right), $$
where $A=\frac{\sqrt{C_3}}{1+C_2^2}, B=C_2,C=\pm\frac{1}{\sqrt{C_3}},m,k  $ are constants.
\end{remark}
\subsection{Spacelike Surface of Revolution of Elliptic Type}
Let $C$ be a spacelike curve in $\text{span}\{e_1,e_3,e_4\}$ parametrized by  arc-length,
$$z(u)=\left(\rho(u),0,f(u),g(u)\right) ,\ u\in I.$$
The orbit of $C$ under the action of the orthogonal transformations of $\mathbb R_1^4$  leaving the timelike plane $Ozt,$
$$A_T=\left[\begin{matrix}\cos v&-\sin v&0&0\\\sin v&\cos v&0&0\\0&0&1&0\\0&0&0&1\end{matrix}\right],\ v\in\mathbb R,$$
is a  surface $M$ given by
\begin{equation}\label{R2}{\rm{\bf X}}(u,v)=\left(\rho(u)\cos v,\rho(u)\sin v,f(u),g(u)\right),\ v\in\mathbb R. \end{equation}
We have the following coefficients of the first fundamental form of surface:
$$g_{11}=(\rho')^2+(f')^2-(g')^2=1,\ g_{12}=0,\ g_{22}=(\rho(u))^2.$$
It follows that  $M$ is a spacelike surface which is called the {\it spacelike surface of revolution of elliptic type} in $\mathbb R_1^4$. We will assume $f'\ne0,g'\ne0$ and $\rho'\ne0$.\\
The system of equations (\ref{hel}) yields $\mathfrak l _1^{\pm}=(\mathfrak l_1,\mathfrak l_2,\mathfrak l_3,1) $, where
$$\mathfrak l_1=\cos v\frac{g'\rho'\pm f'}{(f')^2+\left(\rho'\right)^2},\ \mathfrak l_2=\sin v\frac{g'\rho'\pm f'}{(f')^2+(\rho')^2} ,$$
$$\mathfrak l_3=\frac{g'}{f'}-\frac{\rho'}{f'}\frac{g'\rho'\pm f'}{(f')^2+(\rho')^2} =\frac{g'f'\pm \rho'}{(f')^2+(\rho')^2}.   $$
We are able to define the coefficients of the second fundamental form of $M$ associated to $\mathfrak l_r^{\pm} $-Gauss maps.
$$b_{11}^{\mathfrak l_1^{\pm} }=\langle \mathfrak l_1^{\pm} ,{\rm{\bf X}}_{uu}\rangle=\rho''\frac{g'\rho'\pm f'}{(f')^2+(\rho')^2} +f''\frac{g'f'\mp \rho'}{(f')^2+(\rho')^2}-g'' , $$
$$b_{12}^{\mathfrak l_1^{\pm} }=0,\  b_{22}^{\mathfrak l_1^{\pm} }=\langle \mathfrak l_1^{\pm},{\rm{\bf X}}_{vv} \rangle =-\rho\frac{g'\rho'\pm \rho'}{(f')^2+(\rho')^2}. $$
Solving the equation $\det \left[(b_{ij}^{\mathfrak l_1^{\pm} })-k(g_{ij})\right]=0$, we obtain the following principal curvatures of surface in terms of $\mathfrak l_1^{\pm} $-Gauss map, respectively
\begin{equation}\label{kl1+}k_1^{\mathfrak l_1^{+} }=\rho''\frac{g'\rho'+ f'}{(f')^2+(\rho')^2} +f''\frac{g'f'- \rho'}{(f')^2+(\rho')^2}-g'', \end{equation}
\begin{equation}\label{kl1-}k_1^{\mathfrak l_1^{-} }= \rho''\frac{g'\rho'- f'}{(f')^2+(\rho')^2} +f''\frac{g'f'+ \rho'}{(f')^2+(\rho')^2}-g'', \end{equation}
\begin{equation}\label{kl2+-}k_2^{\mathfrak l_1^{+} }=-\frac{1}{\rho}\frac{g'\rho'+f'}{(f')^2+(\rho')^2};\ \   k_2^{\mathfrak l_1^{-} }=-\frac{1}{\rho}\frac{g'\rho'-f'}{(f')^2+(\rho')^2} .\end{equation}
\begin{theorem}
 If the  surface difined by (\ref{R2}) is totally umbilic, then
$$\rho(u)=\pm\frac{1}{\sqrt C_2}\sinh\left(\sqrt C_2(u-C_1)\right) ,  $$
$$g(u)=\pm\frac{C}{\sqrt C_2}\cosh\left(\sqrt C_2(u-C_1)\right),$$
$$f(u)=\pm\frac{
\sqrt{C+1}}{\sqrt{C_2}}\cosh\left(\sqrt{C_2}(u-C_1)\right)  , $$
where $C,C_1,C_2$ are constants.
\end{theorem}
\begin{proof}
The conditions for $M$ to be totally umbilic,  $$k_1^{\mathfrak l_1^{+} }=k_2^{\mathfrak l_1^{+} }\  \text{and}\  k_1^{\mathfrak l_1^{-} }=k_2^{\mathfrak l_1^{-} },$$ are equivalent to the following system of equations
\begin{equation}\label{helu}\left\{ \begin{aligned}&(f')^2+(\rho')^2-(g')^2=1,&\ &(a)&\\&\rho g''=\rho'g',&\ &(b)&\\&\rho(\rho''f'-f''\rho')=-f'.&\ &(c)&\end{aligned} \right.\end{equation}
It follows from $(b)$ that
\begin{equation}\label{g'p}g'=C\rho,\  C=const.\end{equation}
Differentiating both sides of $(a)$, we obtain
\begin{equation}\label{f'p} f'f''=g'g''-\rho'\rho''.\end{equation}
After multiplying both sides of $(c)$ by $f'$, we have
\begin{equation}\label{rec}\rho(\rho''(f')^2-f'f''\rho')=-(f')^2.\end{equation}
Substitute both (\ref{g'p}) and (\ref{f'p}) to (\ref{rec}), we get a diffirential equation in terms of the function $\rho$
$$\rho(u).\rho''(u)-(\rho'(u))^2+1=0.$$
The claim of the theorem is deduced from Remark \ref{rm1}.
\end{proof}
\begin{remark}
By changing the parameters
$$t=\sqrt{C_2}(u-C_1) ,\ v=v,$$
we get the new parametrization of the totally umbilical spacelike surface of revolution defined by (\ref{R2})
$${\rm{\bf X}}(t,v)=\left(B\sinh t\cos v,B\sinh t\sin v,A\cosh t,CB\cosh t\right), $$
where $A=\pm\frac{\sqrt{C+1}}{\sqrt{C_2}} , B=\pm\frac{1}{\sqrt{C_2}}$ are constants.
\end{remark}
\begin{theorem} If the  surface defined by (\ref{R2}) is  maximal, then
$$\rho(u)=\pm\sqrt{(u-C_1)^2-C},$$
$$g(u)=\frac{C_2}{\sqrt C}\text{\rm arccosh} \left(\frac{u-C_1}{\sqrt C} \right)  +m ,  $$
$$f(u)=\frac{\sqrt{C_2^2-C}}{\sqrt C}\text{\rm arccosh} \left(\frac{u-C_1}{\sqrt C} \right) +k, $$
where $C,C_1,C_2,m,k$ are constants.
\end{theorem}
\begin{proof}
It follows from (\ref{kl1+}), (\ref{kl1-}) and (\ref{kl2+-}) that the conditions for the maximality of $M$ are described by the following equations $$k_1^{\mathfrak l_1^{+} }=-k_2^{\mathfrak l_1^{+} }\  \text{and}\ k_1^{\mathfrak l_1^{-} }=-k_2^{\mathfrak l_1^{-} }.$$
Therefore,  we have a system of equations defining condition of maximal surfaces
\begin{equation}\label{helm}\left\{ \begin{aligned}&(f')^2+(\rho')^2-(g')^2=1,&&(a')&\\&\rho g''=-\rho'g',&&(b')&\\&\rho(\rho''f'-f''\rho')=f'.&&(c')&\end{aligned} \right.\end{equation}
By using the similar method for the solving the system of equations (\ref{helu}) we get an equation in terms of the function $\rho$
$$\rho(u).\rho''(u)+(\rho'(u))^2-1=0.$$
Remark \ref{rm1} yields the result of the theorem.
\end{proof}
\begin{remark}
By changing the parameters
$$t=\text{\rm arccosh}\left(\frac{u-C_1}{\sqrt C} \right)  ,\ v=v,$$
we get a new parametrization of the maximal spacelike surface of revolution of elliptic type defined by (\ref{R2})
$${\rm{\bf X}}(t,v)=\left(A\sinh t\cos v,A\sinh t\sin v,Bt+m,Dt+k \right),$$
where $A=\pm\frac{1}{\sqrt C}, B=\pm\frac{C_2}{\sqrt C},D=\frac{\sqrt{C_2^2-C}}{\sqrt C}$ are constants.
\end{remark}

{\bf ACKNOWLEDGEMENTS.} The author would like to thank Prof. Doan The Hieu  for his valuable guidance and comments.

\address{ Dang Van Cuong\\
Department of Natural Sciences\\
Duy Tan University \\
Danang,  Vietnam}
{dvcuong@duytan.edu.vn}
\end{document}